\documentclass{article}
\usepackage{amsfonts}

\usepackage{mathtools}
\usepackage{amsmath}
\usepackage{amsthm}
\usepackage{amssymb}
\usepackage{mathrsfs}
\usepackage{graphicx}
\usepackage{upgreek}
\def\thefootnote{}
\def\*{\raisebox{.5mm}{*}}

\def\div{{\,\rm div\,}}

\def\R{\mathbb{R}}

\def\H{{\cal H}}

\def\Ga{\Gamma}

\def\<{\langle}
\def\>{\rangle}
\def\om{\omega}
\def\Om{\Omega}

\def\bma{\left[\begin{array}}
\def\ema{\end{array}\right]}
\def\bda{\left|\begin{array}}
\def\eda{\end{array}\right|}
\def\be{\begin{equation}}
\def\ee{\end{equation}}

\newtheorem{thm}{{}\hskip\parindent Theorem}[section]
\newtheorem{lem}{{}\hskip\parindent Lemma}[section]
\newtheorem{pro}{{}\hskip\parindent Proposition}[section]
\newtheorem{exl}{{}\hskip\parindent Example}[section]

\newtheorem{rem}{{}\hskip\parindent Remark}[section]

\large
\title{Stabilization of the critical and subcritical  semilinear inhomogeneous and anisotropic  elastic wave equation } 
\date{}
\author{
Zhen-Hu Ning
\footnotemark[2]
 \ \ Fengyan Yang\thanks{Corresponding author, E-mail address: yangfy16@bjfu.edu.cn}
\footnotemark[3]
\ \ and \ \ Jiacheng Wang
\footnotemark[4]
}
\begin{document}
\maketitle
\renewcommand{\thefootnote}{\fnsymbol{footnote}}
\footnotetext[2]{\scriptsize
Faculty of Information Technology, Beijing University of Technology, Beijing, 100124, China; E-mail address: nzh41034@163.com.}
\footnotetext[3]{\scriptsize School of Sciences Beijing Forestry University, Beijing, 100083, China.
}
\footnotetext[4]{\scriptsize
School of Mathematical Sciences, University of Chinese Academy of Sciences, Beijing 100049, China; e-mail: wjiacheng@amss.ac.cn.

This work is supported by is supported by the Fundamental Research Funds for the Central Universities, NO.BLX201924, the National Science Foundation of China, grant NO.61573342, and Key Research Program of Frontier Sciences, CAS, NO.QYZDJ-SSWSYS011.}
\begin{quote}
\begin{small}
{\bf Abstract} \,\,We prove exponential decay of  the critical and subcritical   semilinear inhomogeneous and anisotropic  elastic wave equation with locally distributed damping on bounded domain. One novelty compared to previous results, is to give a checkable condition of
the  inhomogeneous and anisotropic medias. Another novelty is to establish a framework to study the stability of  the damped  semilinear inhomogeneous and anisotropic elastic wave equation, which  is hard to apply  Carleman estimates  to deal with.  We develop the  Morawetz estimates and the compactness-uniqueness arguments for the  semiliear elastic wave equation  to prove the unique continuation, observability inequality  and  stabilization result.

It is pointing that our proof is  different from the classical method  (See Dehman et al.\cite{ZYY11}, Joly et al.\cite{ZYY16} and Zuazua \cite{ZYY43}), which succeeds for  the subcritical semilinear wave   equation and fails for the critical semilinear wave   equation.
\\[3mm]
\\[3mm]
\end{small}
\end{quote}

\def\theequation{1.\arabic{equation}}
\setcounter{equation}{0}
\section{Some Notations}
\vskip .2cm

 \quad  \ \ Let $O$ be the origin of $\R^n(n\ge3)$  and
 \be r(x)= |x|,\quad  x\in\R^n\ee
be  the standard distance function of $\R^n$.
Moreover, let $\<\cdot,\cdot\>$, $\div$, $\nabla$, $\Delta$ and $I_n=(\delta_{i,j})_{n\times n}$ be the standard inner product of $\R^n$,  the standard divergence operator of $\R^n$, the standard gradient operator of $\R^n$, the standard Laplace operator of $\R^n$ and the unit matrix.

Let $\left(a_{ijkl}\right)_{n\times n\times n\times n}(x)$ be a smooth  tensor function defined on $\R^n$ satisfying
\be a_{ijkl}(x) = a_{jikl}(x) = a_{klij}(x),\quad \ee
for any $x\in \R^n$ and any $1\leq i,j,k,l\leq n$, and  the ellipticity
condition
\be\label{dcds.1} \alpha \sum_{ i,j=1}^n \varepsilon_{ij}  \varepsilon_{ij} \leq \sum_{ i,j,k,l=1}^n a_{ijkl}(x) \varepsilon_{ij}  \varepsilon_{kl}\leq  \beta \sum_{ i,j=1}^n \varepsilon_{ij}  \varepsilon_{ij},\quad x\in \R^n, \ee
for every symmetric tensor $(\varepsilon_{ij})_{n\times n}$, where $\alpha,\beta $ are positive constants.

Let $ \boldsymbol u(x,t)=(u_1(x,t),...,u_n(x,t)): \R^n\times (0,+\infty)\mapsto \R^n$ be a function. Denote
\be u_{i,j}=\frac{\partial u_i}{\partial x_j},\quad  u_{i,t}=\frac{\partial u_i}{\partial t},\quad  u_{i,tt}=\frac{\partial^2 u_i}{\partial t^2},\ee
and
\be\varepsilon_{ij}(\boldsymbol u)=\frac{1}{2}(u_{i,j}+ u_{j,i}),\quad \varepsilon_{ij,k}(\boldsymbol u)=\frac{\partial \varepsilon_{ij}(\boldsymbol u)}{\partial x_k},\ee
for any $x\in \R^n$ and any $1\leq i,j,k\leq n$.  For any $x\in \R^n$ and any $1\leq i,j,k,l\leq n$, we define
\be \sigma_{ij}(\boldsymbol u)=\sum_{k,l=1}^n a_{ijkl}(x)\varepsilon_{kl}(\boldsymbol u).\ee

Denote \be \nabla \boldsymbol u=\left(\nabla u_1,...,\nabla u_n \right)=(u_{j,i})_{n\times n},\ee
 \be |\nabla \boldsymbol u|^2=\sum_{i=1}^n |\nabla  u_i|^2,\ee
\be \boldsymbol\sigma (\boldsymbol u)=\left(\boldsymbol\sigma_1(\boldsymbol u),...,\boldsymbol\sigma_n(\boldsymbol u) \right) =  (\sigma_{ij}(\boldsymbol u))_{n\times n}, \ee
and \be \boldsymbol\varepsilon (\boldsymbol u)=  (\varepsilon_{ij}(\boldsymbol u))_{n\times n}. \ee

Denote
 \be B(h)=\left\{x \Big| \  |x|\leq h\right\},\quad \forall h> 0.\ee
Let $S(r)$ be the sphere in $\R^n$ with radius $r$. Then
\be\left\<X,\frac{\partial}{\partial r}\right\>=0,\quad  \textmd{for }\ \ X\in S(r)_x, x\in \R^n \backslash O.\ee

\def\theequation{2.\arabic{equation}}
\setcounter{equation}{0}
\section{Introduction}
\vskip .2cm

\quad \ \  Let $\Om\subset\R^n$ be a bounded domain with  smooth compact boundary $\Gamma$ and let $\boldsymbol\nu(x)=(\nu_1(x),...,\nu_n(x))$ be the   unit normal vector outside  $\Om$   for  $ x\in  \Ga$.

It is assumed that
$\Ga=\Ga_0\cup\Ga_1$, where $\Ga_0, \Ga_1\subset \Ga$, $\overline\Ga_0\cap\overline\Ga_1=\emptyset$,
$\Ga_0\not=\emptyset$ and
\be \frac{ \partial r }{\partial \boldsymbol\nu }\leq 0,\quad x\in \Ga_0\quad and \quad  \frac{ \partial r }{\partial \boldsymbol\nu }\ge  0,\quad x\in \Ga_1. \ee

Let $\omega$ be an open subset of $\Om$ such that
  \be \omega \supset \bigcup_{x\in  \Ga_1 } \left\{y\in \Om \Big| \  |y-x|<\xi\right\},\ee
  for some $\xi>0$.
\begin{exl}  Let $R_0,R_1,\varepsilon_0$ be  positive constants such that $R_1>R_0,\varepsilon_0<R_1-R_0$. An example can be given by  $\Om= B(R_1) \backslash B(R_0), \om=B(R_1) \backslash B(R_0+\varepsilon_0), \Ga_0=S(R_0)$ and $ \Ga_1=S(R_1)$. \end{exl}

We consider the following system.
\begin{equation} \label{wg.1}
\begin{cases}\boldsymbol u_{tt}-\div\boldsymbol \sigma(\boldsymbol u)+a(x)\boldsymbol u_{t}+\boldsymbol f(\boldsymbol u) =0\qquad (x,t)\in \Om\times
(0,+\infty),\cr
 \boldsymbol u=0\qquad (x,t)\in \Ga_0\times (0,+\infty),\cr
 \boldsymbol\sigma(\boldsymbol u) \boldsymbol\nu^{T}=0 \qquad (x,t)\in \Ga_1\times (0,+\infty),\cr
\boldsymbol u(x,0)=\boldsymbol u_{0}(x),\quad \boldsymbol u_{t}(x,0)=\boldsymbol u_{1}(x)\qquad x\in \Om ,
\end{cases}
\end{equation}
where  $a(x)\in C^{1}(\R^n)$ is a nonnegative function and
\be \label{elastic-f.1} \boldsymbol f(\boldsymbol u)=(|u_1|^{p_1-1}u_1,...,|u_n|^{p_n-1}u_n),  \ee
where
\be 1<p_i\leq \frac{n+2}{n-2}\quad for \quad 1\leq i\leq n.\ee

The energy of the system (\ref{wg.1}) is defined by

\begin{equation}
\label{wg.2}  E(t)=\frac12\int_{\Om}\left(|{\boldsymbol u_t}|^2+\boldsymbol \sigma(\boldsymbol u) \odot \boldsymbol \varepsilon(\boldsymbol u) +2F(\boldsymbol u) \right)dx,
\end{equation}
where
\be F(\boldsymbol u) =\sum_{i=1}^n\frac{1}{p_i+1}|u_{i}|^{p_i+1},\ee
 and $\odot$ is defined by
\be \boldsymbol B\odot \boldsymbol D= \sum_{i,j=1}^n b_{ij}d_{ji}.\ee
for real matrixes $\boldsymbol B=(b_{ij})_{n\times n}$ and $\boldsymbol D=(d_{ij})_{n\times n}$ .
\begin{rem}
The system  (\ref{wg.1}) can be rewritten as
for $1\leq i\leq n$, 
\begin{equation}
\begin{cases} u_{i,tt}-\sum\limits^{n}_{j=1}\sigma_{ij,j}(\boldsymbol u)+a(x)u_{i,t}+|u_i|^{p_i-1}u_i =0\qquad (x,t)\in \Om\times
(0,+\infty),\cr
 u_i\big|_{\Ga_0}=0\qquad t\in(0,+\infty),\cr
 \sum_{j=1}^n \sigma_{ij}(u)\nu_j\big|_{\Ga_1}=0\qquad t\in(0,+\infty),\cr
u_i(x,0)=u_{0i}(x),\quad u_{i,t}(x,0)=u_{1i}(x)\qquad x\in \Om ,\end{cases}
\end{equation}
Then $E(t)$ can be rewritten as
\begin{equation}
E(t)=\frac12\int_{\Om}\left(\sum^{n}_{i=1}u_{i,t}^2+\sum_{ i,j=1}^n\sigma_{ij}(u) \varepsilon_{ij}(u)  +\sum_{i=1}^n\frac{2}{p_i+1}|u_{i}|^{p_i+1}\right)dx.
\end{equation}

\end{rem}

There are a wealth of literatures on the controllability and stabilization of the elastic wave equation. For homogeneous isotropic elastic wave equation, see\cite{A17,A1,A9,A12}.  For homogeneous nonisotropic elastic wave equation, see\cite{A10,A20,A22,A8,A19,A11,B1,A7,A6,A18,B2}. For the inhomogeneous elastic wave equation, see\cite{C1,C2, nzh.10,A23,A21}.

 There exist few literature on the stabilization of the semilinear elastic wave equation. Stabilization of the subcritical semilinear wave equation has been fully studied. See \cite{ZYY6,ZYY9,ZYY11,ZYY16,ZYY19,ZYY43,B2}.  Microlocal analysis given by  \cite{A1,ger1} are the main methods to deal with the stabilization of the semilinear wave equation. However, microlocal analysis doesn't work for the   critical  semilinear wave equation.
As is known, the Morawetz estimate is a simple and effective tool to deal with the energy  estimate on hyperbolic PDEs. See  \cite{w11,w16,w6,w7, w12,w100, B2}. Therefore, we develop the  Morawetz estimates and the compactness-uniqueness arguments  to try to  prove the stabilization of the critical and subcritical semilinear  inhomogeneous  elastic wave equation.

It is pointing that the (elastic) wave equation with Dirichlet/Neumann boundary condition has  special physical meaning, see \cite{DN1,DN2,DN3,DN4}.

\vskip .2cm

The following  assumption is the main assumption.

{\bf Assumption (A)}\,\,\, There exists constant  $\delta>0$ such that for any   $x\in \Om$ and   for every symmetric tensor $(\varepsilon_{ij})_{n\times n}$,
\be \label{2wg.7.4} \sum_{ijkl=1}^n\left((1-\delta)a_{ijkl}-\frac{ r}{2}\frac{\partial a_{ijkl}}{\partial r}\right)\varepsilon_{ij}\varepsilon_{kl}\ge 0.\ee

\begin{rem}
We don't know whether the condition (\ref{2wg.7.4}) is necessary. However from a view of  inhomogeneous and anisotropic wave equation:
\begin{equation}
\begin{cases}u_{tt}-\div A(x)\nabla u=0\qquad (x,t)\in \Om\times (0,+\infty),\cr
u(x,0)=u_0(x), u_t(x,0)=u_1(x)\qquad x\in \Om,\end{cases}
\end{equation}
the condition (\ref{2wg.7.4}) may be a general condition.

Similar with the condition (\ref{2wg.7.4}),  we give the following condition for the inhomogeneous and anisotropic wave equation.
There exists a constant  $ \delta>0$  such that
\be \label{2wg.7.4.1} \left\< \left((1-\delta)A(x)-\frac{ r(x)}{2}\frac{\partial A(x)}{\partial r}\right)X,X \right\>\ge 0\quad for \quad X\in \R^3_x,\ \ x\in \overline{\Om}.\ee

Let $G(x)=A^{-1}(x)$. Let $x\in \R^n,  X,Y\in \R^n_x$ and $Y=G(x)X $.  We deduce that
\begin{eqnarray} &&
 Y^{T}\left(\lambda A(x)-\frac{r}{2}\frac{\partial A(x)}{\partial r}\right)Y \nonumber\\
=&&\left\<G(x)\left(\lambda A(x)-\frac{ r}{2}\frac{\partial A(x)}{\partial r}\right)G(x)X,X\right\> \nonumber\\
=&&\left\<\left(\lambda G(x)+\frac{ r}{2}\frac{\partial (G(x))}{\partial r}\right)X,X\right\> ,
\end{eqnarray}
where $\lambda$ is a contant. It follows from Lemma 3.3 and Lemma 3.4 in \cite{ZHN} that  the condition (\ref{2wg.7.4.1}) is almost equivalent to GCC (geometric control condition).

\end{rem}

\begin{exl}
Let \be a_{ijkl}(x)=\lambda(x)\delta_{ij}\delta_{kl}+\mu(x)(\delta_{ik}\delta_{jl}+\delta_{il}\delta_{jk}), \quad 1\leq i,j,k,l\leq n, \quad x\in \R^n,\ee
where $\lambda(x),\mu(x)\in C^{\infty} (\R^n)$ satisfy
   \be 0<\alpha \leq \mu(x)\leq \beta \quad  and \quad 0<\alpha\leq \lambda(x)+2\mu(x)\leq  \beta,\quad x\in \R^n.\ee

Assume that there exists constant $0<\delta\leq1$ such that
\be (1-\delta)\mu(x)-\frac{r}{2}\frac{\partial \mu(x)}{\partial r}\ge 0,\quad x\in \Om,\ee
and
\be (1-\delta)(\lambda(x)+2\mu(x))-\frac{r}{2}\frac{\partial \left(\lambda(x)+2\mu(x)\right)}{\partial r}\ge 0,\quad x\in \Om.\ee
Then for any $x\in\Om$,
    \begin{eqnarray}&&\ \ \sum_{ijkl=1}^n\left((1-\delta)a_{ijkl}(x)-\frac{r}{2}\frac{\partial a_{ijkl}}{\partial r}\right)\varepsilon_{ij}\varepsilon_{kl}\nonumber\\
&& =\left((1-\delta)\lambda(x)-\frac{r}{2}\frac{\partial \lambda(x)}{\partial r}\right)\left(\sum_{i=1}^n\varepsilon_{ii}\right)^2\nonumber\\
&&\quad+2\left((1-\delta)\mu(x) -\frac{r}{2}\frac{\partial \mu(x)}{\partial r}\right)\sum_{i,j=1}^n\varepsilon_{ij}\varepsilon_{ij}\nonumber\\
&& \ge 0.
\end{eqnarray}
\end{exl}
Well-posedness of the subcritical semilinear  wave equation has been  studied by \cite{BW,GV2,GV8,GV10,J,P1}  and well-posedness of the critical semilinear wave equation has been  studied by \cite{BS,GV10,Gr1,Gr2,Ka1,Ra,SS1,SS2,Str}.
There exists similar results for the nonlinear elastic wave equation. See \cite{Ag1,Pe1,Zha1}.
It is pointing that well-posedness of the critical semilinear wave equation on bounded domain with Dirichlet boundary condition or Neumann boundary condition has been proved by\cite{bur1,bur2}. However,  well-posedness of the critical semilinear wave equation on Riemannian manifold or with variable coefficients is still an open problem. As far as we know,  the  well-posedness of critical semilinear wave equation on Riemannian manifold or with variable coefficients is so hard that there exists no noteworthy study recently. Since we are mainly interested in  stabilization of the system (\ref{wg.1}), we assume the following condition hold  throughout the paper.

Denote
\be H_{\Ga_0}^1(\Om) =\{w\in H^1(\Om),\quad  w\big|_{\Ga_0}=0\}.\ee
{\bf Assumption (S)}\,\,\ Let $E_0>0$ be a constant.  For any $E(0)\leq E_0$, there exists a unique solution of the system (\ref{wg.1}) such that
\be \label{well-posed.1} (u,u_t)\in C\Bigg([0,+\infty),\Big(H_{\Ga_0}^1(\Om)\Big)^n\times \left(L^2(\Om)\right)^n\Bigg).\ee

\begin{rem} If $E_0$ is sufficiently small, the above condition is equivalent to the global existence of the system (\ref{wg.1}) with small initial data.\end{rem}

\begin{thm}\label{t1.2} Let Assumption  ${\bf (A)}$ hold true.
Then there exists  positive constants $C_1,C_2$, which are dependent on $E_0$ given by (\ref{well-posed.1}), such that
 \begin{equation}
\label{ex.1} E(t)\leq C_1 e^{-C_2t}E(0),\quad \forall t>0.
\end{equation}
\end{thm}

 \vskip .5cm
\def\theequation{3.\arabic{equation}}
\setcounter{equation}{0}
\section{Key Lemmas }
\vskip .2cm

\begin{lem}\label{lem.2}
Suppose that $u(x,t)$ solves the system (\ref{wg.1}). Let  $\boldsymbol H= \phi(x)x=\phi(x) \sum_{i=1}^n x_i\frac{\partial}{\partial x_i} =\phi(x)  r\frac{\partial}{\partial r}$,  where  $\phi\in C^{1}(\R^n)$ is a  nonnegative  function. Then
\begin{eqnarray} \label{wg.14.1}
 &&\int_0^T\int_{\Ga}\left( \boldsymbol H(\boldsymbol u) \boldsymbol\sigma(\boldsymbol u) \boldsymbol \nu^{T}\right) d\Ga dt+\frac12\int_0^T\int_{\Ga}
\left(|{\boldsymbol u_t}|^2-\boldsymbol \sigma(\boldsymbol u) \odot \boldsymbol \varepsilon(\boldsymbol u)-2F(\boldsymbol u)\right)\boldsymbol H\cdot \boldsymbol \nu d\Ga dt\nonumber\\
&& \geq\int_{\Omega}\boldsymbol u_{t} \left(\boldsymbol H(\boldsymbol u)\right)^T dx\Big |^T_0+\delta \int_0^T\int_{\Om}\phi(x) \boldsymbol \sigma(\boldsymbol u) \odot \boldsymbol \varepsilon(\boldsymbol u)  dxdt \nonumber\\
&&\quad-C\int_0^T\int_{\Om}r|\nabla \phi| |\nabla \boldsymbol u|^2   dx dt +\int_0^T\int_{\Om} a(x)\boldsymbol u_{t} \left(\boldsymbol H(\boldsymbol u)\right)^Tdxdt\nonumber\\
&&\quad+\frac12\int_0^T\int_{\Om}\left(|\boldsymbol u_t|^2-\boldsymbol \sigma(\boldsymbol u) \odot \boldsymbol \varepsilon(\boldsymbol u)-2F(\boldsymbol u)\right)\div \boldsymbol H dx dt,
\end{eqnarray}

Moreover, assume that $P\in C^1(\R^n): \R^n \mapsto \R^1 $ is a real function. Then
\begin{eqnarray}
\label{wg.14.2}
&&\int_0^T\int_{\Om} \left(|\boldsymbol u_t|^2-\boldsymbol \sigma( \boldsymbol u) \odot \boldsymbol \varepsilon(\boldsymbol u)-\sum_{i=1}^n|u_{i}|^{p_i+1}\right)P dx dt \nonumber\\
\ge &&\int_{\Omega} P \boldsymbol u_{t}\boldsymbol u^T dx\Big |^T_0-C\int_0^T\int_{\Om}  |\nabla P||\boldsymbol u||\nabla\boldsymbol u|dxdt -\int_0^T\int_{\Ga}P \boldsymbol u\boldsymbol\sigma(\boldsymbol u) \boldsymbol \nu^T d\Ga dt\nonumber\\
&&+\frac12\int_{\Om}  P a(x) |\boldsymbol u|^2 dx\Big |^T_0,\end{eqnarray}
where $C$ depends on $\alpha,\beta$, given by (\ref{dcds.1}).

\end{lem}

{\bf Proof}.
 First, we multiply the elastic wave equations in (\ref{wg.1}) by $\left(\boldsymbol H(\boldsymbol u)\right)^T$ and integrate over $\Omega\times
(0,T) $. Note that
\begin{eqnarray}
\boldsymbol \sigma(\boldsymbol u)\odot \left( \nabla\left( \boldsymbol H(\boldsymbol u)\right) \right)&& =\sum_{i,j,m=1}^n\sigma_{ij}(\boldsymbol u) \left(\phi(x) x_m u_{i,m}\right)_j\nonumber\\ &&=\phi(x) \sum_{i,j,m=1}^n\sigma_{ij}(\boldsymbol u) \left( x_m u_{i,m}\right)_j + \sum_{i,j,m=1}^n\sigma_{ij}(\boldsymbol u) \frac{\partial \phi  }{\partial x_j} x_mu_{i,m}\nonumber\\ &&=\phi(x)\left( \sum_{i,j=1}^n\sigma_{ij}(\boldsymbol u) u_{i,j} +\sum_{i,j,m=1}^n\sigma_{ij}(\boldsymbol u) x_mu_{i,jm} \right)\nonumber\\
&&   \quad + \sum_{i,j,m=1}^n\sigma_{ij}(\boldsymbol u) \frac{\partial \phi  }{\partial x_j} x_mu_{i,m}
\nonumber\\ &&=\phi(x)\left( \sum_{i,j=1}^n\sigma_{ij}(\boldsymbol u) \varepsilon_{ij}(\boldsymbol u) +\sum_{i,j,k,l,m=1}^n a_{ijkl}(x)\varepsilon_{kl}(\boldsymbol u) x_m \varepsilon_{ij,m}(\boldsymbol u) \right)\nonumber\\
&&   \quad + \sum_{i,j,m=1}^n\sigma_{ij}(\boldsymbol u) \frac{\partial \phi  }{\partial x_j} x_mu_{i,m}
\nonumber\\ &&=\phi(x)\left( \sum_{i,j=1}^n\sigma_{ij}(\boldsymbol u) \varepsilon_{ij}(\boldsymbol u) +\frac{ r}{2} \frac{\partial}{\partial r}\left(\sum_{i,j=1}^n\sigma_{ij}(\boldsymbol u) \varepsilon_{ij}(\boldsymbol u)\right) -\frac{ r}{2} \sum_{ijkl=1}^n \frac{\partial a_{ijkl}(x) }{\partial r} \varepsilon_{ij}(\boldsymbol u)\varepsilon_{kl}(\boldsymbol u) \right)\nonumber\\
&&   \quad + \sum_{i,j,m=1}^n\sigma_{ij}(\boldsymbol u) \frac{\partial \phi  }{\partial x_j} x_mu_{i,m}
\nonumber\\ &&=\phi(x)  \sum_{ijkl=1}^n\left(a_{ijkl}-\frac{ r}{2}\frac{\partial a_{ijkl}}{\partial r}\right) \varepsilon_{ij}(\boldsymbol u)\varepsilon_{kl}(\boldsymbol u)  +\frac{1}{2}\boldsymbol  H \left(\sum_{i,j=1}^n\sigma_{ij}(\boldsymbol u) \varepsilon_{ij}(\boldsymbol u)\right)  \nonumber\\
&&   \quad  + \sum_{i,j,m=1}^n\sigma_{ij}(\boldsymbol u) \frac{\partial \phi  }{\partial x_j} x_mu_{i,m}
\end{eqnarray}
Hence
\begin{eqnarray}
\boldsymbol \sigma(\boldsymbol u)\odot \left( \nabla\left( \boldsymbol H(\boldsymbol u)\right) \right)&&\ge  \delta \phi(x) \boldsymbol \sigma(\boldsymbol u) \odot \boldsymbol \varepsilon(\boldsymbol u)       -Cr|\nabla \phi| |\nabla \boldsymbol u|^2
\nonumber\\
&&   \quad  +\frac{1}{2}\boldsymbol  H \left( \boldsymbol \sigma(\boldsymbol u) \odot \boldsymbol \varepsilon(\boldsymbol u)\right) \nonumber\\ &&=  \delta \phi(x) \boldsymbol \sigma(\boldsymbol u) \odot \boldsymbol \varepsilon(\boldsymbol u)       -Cr|\nabla \phi| |\nabla \boldsymbol u|^2
\nonumber\\
&&   \quad +\frac12\div \left(\boldsymbol \sigma(\boldsymbol u) \odot \boldsymbol \varepsilon(\boldsymbol u)  H\right) -\frac12\div \boldsymbol H \left( \boldsymbol \sigma(\boldsymbol u) \odot \boldsymbol \varepsilon(\boldsymbol u) \right),
\end{eqnarray}
Therefore
\begin{eqnarray}
0&&=\left(\boldsymbol u_{tt}-\div\boldsymbol\sigma(\boldsymbol u)+a(x)\boldsymbol u_{t}+\boldsymbol f(\boldsymbol u)\right)\left(\boldsymbol H(\boldsymbol u)\right)^T\nonumber\\
&& = \left(\left(\boldsymbol u_{t}\left(\boldsymbol H(\boldsymbol u)\right)^T \right)_t-\frac{1}{2}\div (\boldsymbol u^2_{t}\boldsymbol  H)+\frac{1}{2}\boldsymbol  u^2_{t}\div \boldsymbol  H \right)\nonumber\\
&&  \quad-\left(\div\left(\boldsymbol\sigma(\boldsymbol u)\left(\boldsymbol H(\boldsymbol u)\right)^T \right)-\boldsymbol \sigma(\boldsymbol u)\odot \left( \nabla\left( \boldsymbol H(\boldsymbol u)\right) \right)\right) \nonumber\\
&&  \quad+a(x)\boldsymbol g(u_{t}) \left(\boldsymbol H(\boldsymbol u)\right)^T\nonumber\\
&&  \quad +\boldsymbol H(F(\boldsymbol u))
\nonumber\\
&& = \left(\left(\boldsymbol u_{t}\left(\boldsymbol H(\boldsymbol u)\right)^T \right)_t-\frac{1}{2}\div (\boldsymbol u^2_{t}\boldsymbol  H)+\frac{1}{2}\boldsymbol  u^2_{t}\div \boldsymbol  H \right)\nonumber\\
&&  \quad-\div\left(\boldsymbol\sigma(\boldsymbol u)\left(\boldsymbol H(\boldsymbol u)\right)^T\right) +\boldsymbol \sigma(\boldsymbol u)\odot \left( \nabla\left( \boldsymbol H(\boldsymbol u)\right) \right) \nonumber\\
&&  \quad+a(x)\boldsymbol u_{t}\left(\boldsymbol H(\boldsymbol u)\right)^T\nonumber\\
&&  \quad+ \div (F(\boldsymbol u)\boldsymbol H)-F(\boldsymbol u)\div \boldsymbol  H
\nonumber\\
&& \ge \left(\boldsymbol u_{t}\left(\boldsymbol H(\boldsymbol u)\right)^T \right)_t+\delta \phi(x) \boldsymbol \sigma(\boldsymbol u) \odot \boldsymbol \varepsilon(\boldsymbol u) -Cr|\nabla \phi| |\nabla \boldsymbol u|^2
 \nonumber\\
&&  \quad-\div\left(\boldsymbol\sigma(\boldsymbol u)\left(\boldsymbol H(\boldsymbol u)\right)^T\right) \nonumber\\
&&  \quad+\frac12\left(|\boldsymbol u_{t}|^2-\boldsymbol \sigma(\boldsymbol u) \odot \boldsymbol \varepsilon(\boldsymbol u)-2F(\boldsymbol u)  \right)\div \boldsymbol H\nonumber\\
&&  \quad-\frac12 \div \left(|\boldsymbol u_{t}|^2-\boldsymbol \sigma(\boldsymbol u) \odot \boldsymbol \varepsilon(\boldsymbol u)-2F(\boldsymbol u) \right)\boldsymbol H
\nonumber\\
&&  \quad +a(x)\boldsymbol u_{t} \left(\boldsymbol H(\boldsymbol u)\right)^T .
\end{eqnarray}

In addition, we multiply the wave equation in (\ref{wg.1}) by $Pu$,  and integrate over $\Omega\times
(0,T)$. Note that
\begin{eqnarray}
0&&=\left(\boldsymbol u_{tt}-\div\boldsymbol\sigma(\boldsymbol u)+a(x) \boldsymbol u_{t}+\boldsymbol f(\boldsymbol u)\right)P \boldsymbol u^{T}\nonumber\\
&& =\left(\left(P \boldsymbol  u_{t} \boldsymbol u^T\right)_t-P| \boldsymbol u_{t}|^2\right)\nonumber\\
&&  \quad-\left(\div\left(P\boldsymbol\sigma(u)\boldsymbol u^{T}\right)-P\boldsymbol \sigma(\boldsymbol u)\odot \boldsymbol\varepsilon (\boldsymbol u) -\boldsymbol \sigma(u) \odot  \left( (\nabla P) \boldsymbol u\right) \right)\nonumber\\
&&  \quad+\frac12 (P a(x) \boldsymbol |\boldsymbol u|^{2})_t+P\sum_{i=1}^n|u_{i}|^{p_i+1}\nonumber\\
&& =\left(P \boldsymbol u_{t} \boldsymbol u\right)_t-\div\left(\boldsymbol\sigma(u) P \boldsymbol u^{T}\right)+\boldsymbol \sigma(\boldsymbol u) \odot  \left(( \nabla P) \boldsymbol u\right) +\frac12 (P a(x) \boldsymbol |\boldsymbol u|^{2} )_t \nonumber\\
&&  \quad-P\left(|\boldsymbol u_t|^2-\boldsymbol \sigma(\boldsymbol u) \odot \boldsymbol \varepsilon(\boldsymbol u)-\sum_{i=1}^n|u_{i}|^{p_i+1}\right)
\nonumber\\
&& \ge \left(P \boldsymbol u_{t} \boldsymbol u\right)_t-\div\left(\boldsymbol\sigma(u) P \boldsymbol u^{T}\right)-C |\nabla P||\boldsymbol u||\nabla\boldsymbol u| +\frac12 \left(P a(x) \boldsymbol |\boldsymbol u|^{2}\right)_t \nonumber\\
&&  \quad-P\left(|\boldsymbol u_t|^2-\boldsymbol \sigma(\boldsymbol u) \odot \boldsymbol \varepsilon(\boldsymbol u)-\sum_{i=1}^n|u_{i}|^{p_i+1}\right).
\end{eqnarray}
The equality  (\ref{wg.14.2}) follows from Green's formula.$\Box$

\begin{lem}\label{wg-kl.22}Let $u(x,t)$ solve the system (\ref{wg.1}). Then
\be \label{wg-kl.23} E(t) \Big|_0^T=-\int_0^T \int_{\Om }a(x) |\boldsymbol u_{t}|^2dx dt ,\ee
 \end{lem}
which implies $E(t)$ is decreasing.

{\bf Proof}.  Multiplying the equation in (\ref{wg.1}) by $\boldsymbol u_{t}$, and integrating over $\Om\times
(0,T)$, the equality (\ref{wg-kl.23}) follows from Green's formula. $\Box$

 \begin{pro}\label{elastic.2} Let Assumption  ${\bf (A)}$ hold true. Then  there exists $T_0\ge0$ such that for any $T>T_0$, the only solution $(u,u_t)\in C\left([0,T],\left(H^1(\Om)\right)^n\times \left(L^2(\Om)\right)^n\right)$ to the system
\begin{equation}
\label{elastic.3} \begin{cases}\boldsymbol u_{tt}-\div \boldsymbol\sigma(\boldsymbol u) +\boldsymbol f(\boldsymbol u) =0\qquad (x,t)\in \Om \times
(0,T),\cr \boldsymbol u\large|_{\Ga_0}=0\qquad t\in(0,+\infty),\cr
 \boldsymbol\sigma(\boldsymbol u)\boldsymbol\nu\large|_{\Ga_1}=0\qquad t\in(0,+\infty),\cr \boldsymbol u_t=0\qquad(x,t)\in \omega \times
(0,T),\end{cases}
\end{equation}
where $\boldsymbol f(\boldsymbol u)$ is given by (\ref{elastic-f.1}), is the trivial one $\boldsymbol u\equiv 0$. \end{pro}
{\bf Proof } \  \   
Let $a(x)\equiv 0$, it follows from (\ref{wg-kl.23}) that
\be E(t)=E(0),\quad t>0.\ee

Let $\boldsymbol H =x$ and $a(x)\equiv 0$. It follows from (\ref{wg.14.1}) that
\begin{eqnarray} \label{wg.14.10}
 &&\int_0^T\int_{\partial\Om}\left( \boldsymbol H(\boldsymbol u) \boldsymbol\sigma(\boldsymbol u) \boldsymbol \nu^{T}\right) d\Ga dt+\frac12\int_0^T\int_{\partial\Om}
\left(|{\boldsymbol u_t}|^2-\boldsymbol \sigma(\boldsymbol u) \odot \boldsymbol \varepsilon(\boldsymbol u)-2F(\boldsymbol u)\right)\boldsymbol H\cdot \boldsymbol \nu d\Ga dt\nonumber\\
&& \geq\int_{\Omega}\boldsymbol u_{t} \left(\boldsymbol H(\boldsymbol u)\right)^T dx\Big |^T_0+\delta \int_0^T\int_{\Om}\phi(x) \boldsymbol \sigma(\boldsymbol u) \odot \boldsymbol \varepsilon(\boldsymbol u)  dxdt \nonumber\\
&&\quad-C\int_0^T\int_{\Om}r|\nabla \phi| |\nabla \boldsymbol u|^2   dx dt \nonumber\\
&&\quad+\frac12\int_0^T\int_{\Om}\left(|\boldsymbol u_t|^2-\boldsymbol \sigma(\boldsymbol u) \odot \boldsymbol \varepsilon(\boldsymbol u)-2F(\boldsymbol u)\right)\div \boldsymbol H dx dt\nonumber\\
&& =\int_{\Omega}\boldsymbol u_{t} \left(\boldsymbol H(\boldsymbol u)\right)^T dx\Big |^T_0+\delta \int_0^T\int_{\Om}\phi(x) \boldsymbol \sigma(\boldsymbol u) \odot \boldsymbol \varepsilon(\boldsymbol u)  dxdt \nonumber\\
&&\quad-C\int_0^T\int_{\Om}r|\nabla \phi| |\nabla \boldsymbol u|^2   dx dt \nonumber\\
&&\quad+\frac{n}{2}\int_0^T\int_{\Om}\left(|\boldsymbol u_t|^2-\boldsymbol \sigma(\boldsymbol u) \odot \boldsymbol \varepsilon(\boldsymbol u)-\sum_{i=1}^n|u_{i}|^{p_i+1}\right) dx dt\nonumber\\
&&\quad + \int_0^T\int_{\Om}\sum_{i=1}^n\frac{(p_i-1)n}{2(p_i+1)}|u_{i}|^{p_i+1}dxdt.
\end{eqnarray}

Denote
 \be \delta_c=\min_{1\leq i\leq n}\{\delta,\frac{( p_i-1)n}{2( p_i+1)}\}.\ee
Let $P=\frac{n-\delta_c}{2}$ and $a(x)\equiv 0$.
Substituting the formula
(\ref{wg.14.2}) into the formula  (\ref{wg.14.10}), we obtain
\begin{eqnarray}
 \label{elastic-f.3}
 &&\int_0^T\int_{\partial\Om}\left( \boldsymbol H(\boldsymbol u) \boldsymbol\sigma(\boldsymbol u) \boldsymbol \nu^{T}\right) d\Ga dt+\frac12\int_0^T\int_{\partial\Om}
\left(|{\boldsymbol u_t}|^2-\boldsymbol \sigma(\boldsymbol u) \odot \boldsymbol \varepsilon(\boldsymbol u)-2F(\boldsymbol u)\right)\boldsymbol H\cdot \boldsymbol \nu d\Ga dt\nonumber\\
&& \geq \int_{\Om}\boldsymbol u_{t}\left(\boldsymbol H(\boldsymbol u) +P\boldsymbol u\right)^Tdx\Big |^T_0\nonumber\\
 &&\quad+\frac{\delta_c}{2}\int_0^T\int_{\Om}\left(|\boldsymbol u_t|^2+\boldsymbol \sigma(\boldsymbol u) \odot \boldsymbol \varepsilon(\boldsymbol u)+2F(\boldsymbol u)\right) dx dt,
\end{eqnarray}

Note that  $u |_{\Ga_0}=0$, then for $1\leq i,j,m\leq n$, \be u_{i,m}\nu_j=u_{i,\boldsymbol\nu}\nu_m\nu_j=u_{i,j}\nu_m, \quad x\in \Ga_0.\ee
Hence
\begin{eqnarray} \boldsymbol H(\boldsymbol u) \boldsymbol\sigma(\boldsymbol u) \boldsymbol \nu^{T} &&=\sum^n_{i,j,m=1}x_mu_{i,m}\sigma_{ij}(\boldsymbol u) \nu_j\nonumber\\
 &&= \sum^n_{i,j,m=1}u_{i,j}\sigma_{ij}(\boldsymbol u) x_m\nu_m
 \nonumber\\
 &&= \sum^n_{i,j,m=1}\varepsilon_{ij}(\boldsymbol u)\sigma_{ij}(\boldsymbol u) x_m\nu_m  \nonumber\\
 &&=\boldsymbol \sigma(\boldsymbol u) \odot \boldsymbol \varepsilon(\boldsymbol u) (\boldsymbol H\cdot \boldsymbol \nu),\quad x\in \Ga_0.\end{eqnarray}

With
\be \boldsymbol u_t=\boldsymbol \sigma(\boldsymbol u)\boldsymbol \nu=0,\quad x\in \Ga_1, \ee
and
\be \frac{ \partial r }{\partial\boldsymbol \nu }\leq 0,\quad x\in \Ga_0\quad and \quad  \frac{ \partial r }{\partial \boldsymbol \nu }\ge  0,\quad x\in \Ga_1,  \ee
we obtain
\begin{eqnarray} \label{elastic-f.2}
 &&\quad \int_0^T\int_{\Ga}\left( (\boldsymbol H(\boldsymbol u)+P\boldsymbol u) \boldsymbol\sigma(\boldsymbol u) \boldsymbol \nu^{T}\right) d\Ga dt\nonumber\\
 &&\quad +\frac12\int_0^T\int_{\Ga}
\left(|{\boldsymbol u_t}|^2-\boldsymbol \sigma(\boldsymbol u) \odot \boldsymbol \varepsilon(\boldsymbol u)-2F(\boldsymbol u)\right)\boldsymbol H\cdot \boldsymbol \nu d\Ga dt\nonumber\\
&&=\quad \int_0^T\int_{\Ga_0}\left( (\boldsymbol H(\boldsymbol u)+P\boldsymbol u) \boldsymbol\sigma(\boldsymbol u) \boldsymbol \nu^{T}\right) d\Ga dt\nonumber\\
 &&\quad +\frac12\int_0^T\int_{\Ga_0}
\left(|{\boldsymbol u_t}|^2-\boldsymbol \sigma(\boldsymbol u) \odot \boldsymbol \varepsilon(\boldsymbol u)-2F(\boldsymbol u)\right)\boldsymbol H\cdot \boldsymbol \nu d\Ga dt\nonumber\\
&&\quad +\int_0^T\int_{\Ga_1}\left( (\boldsymbol H(\boldsymbol u)+P\boldsymbol u) \boldsymbol\sigma(\boldsymbol u) \boldsymbol \nu^{T}\right) d\Ga dt\nonumber\\
 &&\quad +\frac12\int_0^T\int_{\Ga_1}
\left(|{\boldsymbol u_t}|^2-\boldsymbol \sigma(\boldsymbol u) \odot \boldsymbol \varepsilon(\boldsymbol u)-2F(\boldsymbol u)\right)\boldsymbol H\cdot \boldsymbol \nu d\Ga dt\nonumber\\
&&=\quad \frac12\int_0^T\int_{\Ga_0}
\boldsymbol \sigma(\boldsymbol u) \odot \boldsymbol \varepsilon(\boldsymbol u) (\boldsymbol H\cdot \boldsymbol \nu) d\Ga dt\nonumber\\
&&\quad -\frac12\int_0^T\int_{\Ga_1}
\left(\boldsymbol \sigma(\boldsymbol u) \odot \boldsymbol \varepsilon(\boldsymbol u)+2F(\boldsymbol u)\right)\boldsymbol H\cdot \boldsymbol \nu d\Ga dt\nonumber\\
&&\leq 0.
\end{eqnarray}

It follows from the   Korn's inequality with Dirichlet boundary conditions \cite{Jn1} that
\begin{eqnarray} \label{elastic.5}
&& \int_{\Om}
|\nabla \boldsymbol u|^2 dx  \nonumber\\
&& \leq  C\int_{\Om}
 \boldsymbol \varepsilon(\boldsymbol u)  \boldsymbol \varepsilon(\boldsymbol u) dx\nonumber\\
&& \leq C\int_{\Om}
\boldsymbol \sigma(u)\boldsymbol \varepsilon(\boldsymbol u) dx.
\end{eqnarray}
Substituting (\ref{elastic-f.2})  into (\ref{elastic-f.3}), we obtain
\begin{eqnarray}
 &&\int_0^T E(t)dt \leq CE(0),
\end{eqnarray}
which implies
\be (T-C)E(0)\leq 0.\ee
The assertion (\ref{elastic.3}) holds true.$\Box$

By a similar proof with Proposition \ref{elastic.2}, the following proposition holds.
\begin{pro}\label{elastic-f.4} Let Assumption  ${\bf (A)}$ hold true. Then  there exists $T_0\ge0$ such that for any $T>T_0$, the only solution $(\boldsymbol u,\boldsymbol u_t)\in C\left([0,T],\left(H^1(\Om)\right)^n\times \left(L^2(\Om)\right)^n\right)$ to the system
\begin{equation}
\label{elastic.3} \begin{cases}\boldsymbol u_{tt}-\div \boldsymbol \sigma(\boldsymbol u)  =0\qquad (x,t)\in \Om \times
(0,T),\cr \boldsymbol u\big|_{\Ga_0}=0\qquad t\in(0,+\infty),\cr
 \boldsymbol \sigma(\boldsymbol u)\boldsymbol \nu\big|_{\Ga_1}=0\qquad t\in(0,+\infty),\cr \boldsymbol u_t=0\qquad(x,t)\in \omega \times
(0,T),\end{cases}
\end{equation}
 is the trivial one $u\equiv 0$. \end{pro}

 \vskip .5cm
\def\theequation{4.\arabic{equation}}
\setcounter{equation}{0}
\section{Proofs of Theorem \ref{t1.2} }
\vskip .2cm

\begin{lem}\label{TheL.2} \ Let Assumption{\bf (A)} hold true. Let $u$ solve
the system $(\ref{wg.1})$. Then there exists a positive
constant $C$ such that
 \begin{eqnarray}
 \label{The.3}
 &&\int_0^T E(t) dt\leq C \int_0^T\int_{\Om }a(x) |\boldsymbol u_{t}|^2dxdt +C\int_0^T\int_{\Om} |\boldsymbol  u|^2dxdt
 \end{eqnarray}
\end{lem}for sufficiently large $T$.

{\bf Proof}. It follows from classical  Korn's inequality\cite{A26} that
\begin{eqnarray}
\label{elastic.4}
&& \int_{\omega}
|\nabla \boldsymbol u|^2 dx \nonumber\\
&& \leq C\int_{\omega}
\left(|\boldsymbol u|^2+\boldsymbol \varepsilon(\boldsymbol u) \odot \boldsymbol \varepsilon(\boldsymbol u) \right)dx,\end{eqnarray}
and  the   Korn's inequality with Dirichlet boundary conditions \cite{Jn1} that
\begin{eqnarray} \label{elastic.5}
&& \int_{\Om}
|\nabla \boldsymbol u|^2 dx  \nonumber\\
&& \leq  C\int_{\Om}
 \boldsymbol \varepsilon(\boldsymbol u)  \odot\boldsymbol \varepsilon(\boldsymbol u) dx\nonumber\\
&& \leq C\int_{\Om}
\boldsymbol \sigma(\boldsymbol u)\odot\boldsymbol \varepsilon(\boldsymbol u) dx.
\end{eqnarray}

Let $\widehat{\omega} \subset \Om$ be an bounded open set with smooth boundary such that
\be \Ga_1 \subset \overline{\widehat{\omega}}  \quad and \quad \left(\overline{\widehat{\omega}}  \backslash \Ga_1\right) \subset  \omega .\ee
 Let $\phi\in C^{\infty}(\R^n)$ be a  nonnegative  function such that
\be\phi=1,x\in \Om \backslash\omega\quad and \quad \phi=0,x\in \widehat{\omega}. \ee
Let
\be \boldsymbol H=\phi(x)x.\ee
It follows from  (\ref{wg.14.1}) that
\begin{eqnarray} \label{elastic.12}
 &&\int_0^T\int_{\Ga_0}\left( \boldsymbol H(\boldsymbol u) \boldsymbol\sigma(\boldsymbol u) \boldsymbol \nu^{T}\right) d\Ga dt+\frac12\int_0^T\int_{\Ga_0}
\left(|{\boldsymbol u_t}|^2-\boldsymbol \sigma(\boldsymbol u) \odot \boldsymbol \varepsilon(\boldsymbol u)-2F(\boldsymbol u)\right)\boldsymbol H\cdot \boldsymbol \nu d\Ga dt\nonumber\\
&& \geq\int_{\Omega}\boldsymbol u_{t} \left(\boldsymbol H(\boldsymbol u)\right)^T dx\Big |^T_0+\delta \int_0^T\int_{\Om}\phi(x) \boldsymbol \sigma(\boldsymbol u) \odot \boldsymbol \varepsilon(\boldsymbol u)  dxdt \nonumber\\
&&\quad-C\int_0^T\int_{\Om}r|\nabla \phi| |\nabla \boldsymbol u|^2   dx dt +\int_0^T\int_{\Om} a(x)\boldsymbol u_{t} \left(\boldsymbol H(\boldsymbol u)\right)^Tdxdt\nonumber\\
&&\quad+\frac12\int_0^T\int_{\Om}\left(|\boldsymbol u_t|^2-\boldsymbol \sigma(\boldsymbol u) \odot \boldsymbol \varepsilon(\boldsymbol u)-\sum_{i=1}^n|u_{i}|^{p_i+1}\right)\div \boldsymbol H dx dt\nonumber\\
&&\quad + \int_0^T\int_{\Om}\sum_{i=1}^n\frac{(p_i-1)\div \boldsymbol H }{2(p_i+1)}|u_{i}|^{p_i+1}dxdt .
\end{eqnarray}

Note that
\be \div \boldsymbol H=n,\quad x\in \Om\backslash \omega.\ee
Denote
 \be \delta_c=\min_{1\leq i\leq n}\{\delta,\frac{( p_i-1)n}{2( p_i+1)}\}.\ee
Let $P=\frac{1}{2}\left(\div\boldsymbol  H-\phi \delta_c\right)$, substituting (\ref{wg.14.2})  into (\ref{elastic.2}), with (\ref{elastic.4})  and (\ref{elastic.5}), we obtain
\begin{eqnarray}\label{kl.1}
&&\int_0^T\int_{ \Om\backslash\om}\left(\boldsymbol |\boldsymbol u_{t}|^2+\boldsymbol \sigma(\boldsymbol u)\odot\boldsymbol \varepsilon(\boldsymbol u)+2F(\boldsymbol u)\right) dx dt\nonumber\\
&&\leq  C(E(0)+E(T))+ C\int_0^T\int_{\om}\left(|\boldsymbol u_{t}|^2+\boldsymbol \sigma(\boldsymbol u)\odot\boldsymbol \varepsilon(\boldsymbol u)+\sum_{i=1}^n|u_{i}|^{p_i+1}\right) dx dt
 \nonumber\\
&&\qquad+C\int_0^T\int_{\Om}a(x)\left(|\boldsymbol  u_{t}|^2+\boldsymbol \sigma(\boldsymbol u)\odot\boldsymbol \varepsilon(\boldsymbol u)\right) dx dt \nonumber\\
&&\qquad+\int_0^T\int_{\Om}\left(C_\epsilon |\boldsymbol u|^2+ \epsilon\boldsymbol \sigma(\boldsymbol u)\odot\boldsymbol \varepsilon(\boldsymbol u)\right) dx dt.\end{eqnarray}
Hence
\begin{eqnarray}
&&\int_0^T\int_{ \Om\backslash\om}\left(|\boldsymbol u_{t}|^2+\boldsymbol \sigma(\boldsymbol u)\odot\boldsymbol \varepsilon(\boldsymbol u)+2F(\boldsymbol u)\right) dx dt\nonumber\\
&&\leq  C(E(0)+E(T))+ C\int_0^T\int_{\Om}a(x)\left(|\boldsymbol u_{t}|^2+\boldsymbol \sigma(\boldsymbol u)\odot\boldsymbol \varepsilon(\boldsymbol u)+\sum_{i=1}^n|u_{i}|^{p_i+1}\right) dx dt
\nonumber\\
&&\quad  +C\int_0^T\int_{\Om}|\boldsymbol  u|^2 dx dt.\end{eqnarray}
Therefore
\begin{eqnarray}\label{kl.2}
&&\int_0^T\int_{\Om}\left(\boldsymbol u_{t}^2+\boldsymbol \sigma(\boldsymbol u)\odot\boldsymbol \varepsilon(\boldsymbol u)+2F(\boldsymbol u)\right) dx dt\nonumber\\
&&\leq  C(E(0)+E(T))+ C\int_0^T\int_{\Om}a(x)\left(|\boldsymbol u_{t}|^2+\boldsymbol \sigma(\boldsymbol u)\odot\boldsymbol \varepsilon(\boldsymbol u)+\sum_{i=1}^n|u_{i}|^{p_i+1}\right) dx dt
\nonumber\\
&&\quad +C\int_0^T\int_{\Om} |\boldsymbol u|^2 dx dt.\end{eqnarray}

Let $P=a(x)$ in the equality (\ref{wg.14.2}), we obtain
\begin{eqnarray}
&&\int_0^T\int_{\Om} a(x)\left(\boldsymbol \sigma(\boldsymbol u)\odot\boldsymbol \varepsilon(\boldsymbol u)+2F(\boldsymbol u)\right)dxdt \nonumber\\
&&\leq C(E(0)+E(T))+C\int_0^T\int_{\Om} a(x)\left(|\boldsymbol u_{t}|^2+|\boldsymbol u|^2\right)dxdt\nonumber\\
&& \quad +C \int_0^T\int_{\Om}a(x) |\boldsymbol u_{t}|^2dxdt \nonumber\\
&& \quad +\int_0^T\int_{\Om}\left(C_\epsilon|\boldsymbol u|^2+\epsilon\boldsymbol \sigma(\boldsymbol u)\odot\boldsymbol \varepsilon(\boldsymbol u)\right) dx dt.\end{eqnarray}
With (\ref{kl.2}), we obtain
\begin{eqnarray}\label{kl.3}
\int_0^T E(t) dt &&\leq C(E(0)+E(T))+C \int_0^T\int_{\Om}a(x) |\boldsymbol u_{t}|^2dxdt\nonumber\\
&& \quad +C\int_0^T\int_{\Om}|\boldsymbol  u|^2dxdt.\end{eqnarray}

It follows from (\ref{wg-kl.23}) that
\begin{eqnarray}\label{kl.4}
C(E(0)+E(T)) &&= 2CE(T)+C\int_0^T\int_{\Om} a(x)|\boldsymbol u_{t}|^2dxdt.\end{eqnarray}
Note that $E(t)$ is decreasing, then, for $T\ge 4C$
\begin{eqnarray}\label{kl.5}
&&\int_0^T E(t) dt\leq C \int_0^T\int_{\Om}a(x) |\boldsymbol u_{t}|^2dxdt +C\int_0^T\int_{\Om} |\boldsymbol  u|^2dxdt.\end{eqnarray}$\Box$

\begin{lem}\label{wg-kl.29} Let Assumption  ${\bf (A)}$ hold true. Let $u(x,t)$  solve the system (\ref{wg.1}). Then for any $E(0) \leq E_0$,
\begin{eqnarray}\label{wg-kl.30}
&&\int_0^T E(t) dt\leq C(E_0,T) \int_0^T\int_{\Om}a(x) |\boldsymbol u_{t}|^2dxdt,\end{eqnarray}  for sufficiently large $T$.\end{lem}
{\bf Proof}.
We apply compactness-uniqueness arguments to prove the conclusion. It follows from (\ref{TheL.2}) that
 \begin{eqnarray}  \label{wg-kl.85}
&&\int_0^T E(t) dt\leq C \int_0^T\int_{\Om}a(x) |\boldsymbol u_{t}|^2dxdt +C\int_0^T\int_{\Om} |\boldsymbol  u|^2dxdt.\end{eqnarray}
Then, if the estimate (\ref{wg-kl.30}) doesn't hold true, there exist $\Big\{u_k\Big\}_{k=1}^{\infty}$ such that
 \be \label{scd.08} E_k(0)\leq E_0,\ee
 where
 \be  E_k(t)=\frac12\int_{\Om}\left(|\boldsymbol u_{k,t}|^2+\boldsymbol \sigma(\boldsymbol u_k) \odot\boldsymbol \varepsilon(\boldsymbol u_k)  +2F(\boldsymbol u_{k})\right),\ee
  and
  \be \label{wg-kl.81}\int_0^T\int_{\Om}   |\boldsymbol u_{k}|^2dxdt \ge k \int_0^T\int_{\Om}a(x)  |\boldsymbol u_{k,t}|^2dxdt.\ee
With (\ref{wg-kl.23}), we have
    \be \label{anapde.07} E_k(t)\leq  E_0,\quad 0\leq t\leq T.\ee
   and  \be \int_0^T E_k(t)dt \leq  T E_0.\ee
   Therefore, there exists $\hat{\boldsymbol u}_0$ and a subset of $\Big\{\boldsymbol u_k\Big\}_{k=1}^{\infty}$, still denoted by $\Big\{\boldsymbol u_k\Big\}_{k=1}^{\infty}$,  such that
   \be \label{anapde.01} \boldsymbol u_k\rightarrow \hat{\boldsymbol u}_0 \ \ weakly \ \  in  \ \    \left(H^1\left(\Om\times(0,T)\right)\right)^n ,\ee
   and
    \be \boldsymbol u_k\rightarrow \hat{\boldsymbol u}_0 \ \ strongly \ \ in  \ \  \left( L^2\left(\Om\times(0,T)\right)\right)^n .\ee

{\bf Case a:}
 \be \label{wg-kl.80}\int_0^T\int_{\Om} |\hat{\boldsymbol u}_{0}|^2dxdt>0.\ee

Denote
\be \label{anapde.02} \boldsymbol q_i=\frac{2n}{(n-2) p_i},\quad \boldsymbol q_i^{*}=\frac{ q_i}{ q_i-1}, \ee
for $1\leq i\leq n$.
Since $1< p_i\leq \frac{n+2}{n-2}$, then
\be \frac{2n}{n+2}\leq q_i,q_i^{*}\leq \frac{2n}{n-2},\quad .\ee
Note that\be\frac{1}{q_i}+\frac{1}{q_i^{*}}=1,\quad \ee
Then, $L^{ q_i^{*}}\left(\Om\right)$ is the dual space of $L^{ q_i}\left(\Om\right)$.

Note that
\be H^{1}\left(\Om\right)\hookrightarrow L^{\frac{2n}{n-2}}\left(\Om\right).\ee
Therefore, it follows from (\ref{anapde.07}) that for $1\leq i\leq n$
\be \{ |u_{ki}|^{p_i-1}u_{ki}\}\ are \ bounded \ in\  L^{\infty}([0,T],L^{ q_i}(\Om)).\ee
 Then
 \be \{|u_{ki}|^{p_i-1}u_{ki} \}\ are \ bounded \ in\  L^{ q_i} \left(\Om\times(0,T)\right).\ee
 Hence for $1\leq i\leq n$
  \be \label{scd.24} |u_{ki}|^{p_i-1}u_{ki} \rightarrow |\hat{u}_{0i}|^{p_i-1}\hat{u}_{0i}\ \ weakly \ \  in  \ \ L^{q_i} \left(\Om\times(0,T)\right).\ee

 It follows from (\ref{wg-kl.81}) that
 \be a(x)\hat{\boldsymbol u}_{0t}=0\qquad(x,t)\in \R^n\times
(0,T).\ee
 Therefore, with  (\ref{anapde.01}) and (\ref{scd.24}),  we obtain \begin{equation}
\label{wg-kl.32} \begin{cases}\hat{\boldsymbol u}_{0tt}-\div \boldsymbol \sigma(\hat{\boldsymbol u}_0) +\boldsymbol f(\hat{\boldsymbol u}_0) =0\qquad (x,t)\in \Om \times
(0,T),\cr \hat{\boldsymbol u}_0\large|_{\Ga_0}=0\qquad t\in(0,+\infty),\cr
 \boldsymbol \sigma(\hat{\boldsymbol u}_0)\boldsymbol \nu\large|_{\Ga_1}=0\qquad t\in(0,+\infty),\cr \hat{ \boldsymbol u}_{0t}=0\qquad(x,t)\in \omega \times
(0,T),\end{cases}
\end{equation}
where $f(\hat{u}_0)$ is given by (\ref{elastic-f.1}).
It follows from  Proposition \ref{elastic.2} that
\be \label{wg-kl.82} \hat{\boldsymbol u}_0\equiv0,\qquad (x,t)\in \Om\times
(0,T), \ee
which contradicts (\ref{wg-kl.80}).

{\bf Case b:}
 \be \label{scd.05} \hat{\boldsymbol u}_0\equiv 0 \quad on \quad  \Om\times(0,T). \ee

Denote
\be \label{scd.06}\boldsymbol v_k=  \boldsymbol u_k\Big/ \sqrt{c_k} \quad for\quad  k\ge 1,\ee
where
\be  \label{scd.27}c_k=\int_0^T\int_{\Om}  |\boldsymbol u_{k}|^2 dx dt.\ee
Then $v_k$ satisfies for $1\leq i\leq n$,
\be  \label{wg-kl.90}\boldsymbol v_{k,tt}-\div\boldsymbol \sigma(\boldsymbol v_k)+a (x)\boldsymbol v_{k,t}+
\frac{\boldsymbol f(\boldsymbol u_{k} )}{\sqrt{c_k}}=0\qquad (x,t)\in \Om\times
(0,T),\ee
and
\be\label{wg-kl.89}\int_0^T\int_\Om  |\boldsymbol v_{k}|^2 dx dt=1.\ee

It follows from (\ref{wg-kl.81}) that
\be \label{wg-kl.88}  1\ge k  \int_0^T  \int_{\Om }a(x) \boldsymbol v_{kt}^2dx dt .\ee
Therefore, It follows from  (\ref{wg-kl.85}) that
\begin{eqnarray}  \label{wg-kl.86}
\widehat{E}_k(0)&&\leq 1+\frac{1}{k}\leq 2,\end{eqnarray}
where
 \be\widehat{E}_k(t)= \frac12\int_{\Om}\left(|\boldsymbol v_{k,t}|^2+\boldsymbol \sigma(\boldsymbol v_k)\odot \boldsymbol \varepsilon(\boldsymbol u_k)  +\sum_{i=1}^n\frac{2}{ p_i+1}| u_{ki}|^{ p_i-1}| v_{ki}|^{2}\right).\ee
 Hence, there exists $v_0$ and a subset of $\Big\{\boldsymbol v_k\Big\}_{k=1}^{\infty}$, still denoted by $\Big\{\boldsymbol v_k\Big\}_{k=1}^{\infty}$, such that
   \be \boldsymbol v_k\rightarrow\boldsymbol  v_0 \ \ weakly \ \  in  \ \   \left(H^1\left(\Om\times(0,T)\right)\right)^n,\ee
   and
    \be \boldsymbol v_k\rightarrow \boldsymbol v_0 \ \ strongly \ \ in  \ \ \left(L^2\left(\Om\times(0,T)\right)\right)^n.\ee

It follows from  (\ref{wg-kl.23}), (\ref{scd.06})  and (\ref{wg-kl.86}) that
\be \label{scd.26}\widehat{E}_k(t)\leq \widehat{E}_k(0)\leq 2,\quad \forall 0\leq t\leq T.\ee
 Let $q_i,q_i^{*}$ be given by (\ref{anapde.02}).  Note that
 \be \label{scd.07}H^{1}\left(\Om\right)\hookrightarrow L^{\frac{2n}{n-2}}\left(\Om\right).\ee
Therefore, it follows from (\ref{scd.26}) that for $1\leq i\leq n$
\begin{eqnarray} \int_0^T\int_{\Om}&&  \sum_{i=1}^n  \left(|u_{ki}|^{p_i-1}v_{ki}\right)^{ q_i} dx_gdt\nonumber\\ &&=\int_0^T\int_{\Om} \sum_{i=1}^n c_k^{\frac{ q_i( p_i-1)}{2}} | v_{ki}|^{\frac{2n}{n-2}} dxdt\nonumber\\ && \leq\sum_{i=1}^n c_k^{\frac{ q_i( p_i-1)}{2}} C(T).\end{eqnarray}
With (\ref{scd.05}) and (\ref{scd.27}),  we obtain
  \be \label{anapde.03}\lim_{k\rightarrow +\infty}\int_0^T\int_{\Om} \sum_{i=1}^n  \left(|u_{ki}|^{p_i-1}v_{ki}\right)^{ q_i} dxdt = 0.\ee

 It follows from (\ref{wg-kl.88}) that
 \be a(x)\boldsymbol v_{0t}=0\qquad(x,t)\in \R^n\times
(0,T).\ee
Therefore, It follows from (\ref{wg-kl.90}) and (\ref{anapde.03}) that
\begin{equation}
\ \begin{cases} \boldsymbol v_{0tt}-\div\boldsymbol \sigma(\boldsymbol v_0)=0\qquad (x,t)\in \Om\times(0,T),
\cr \boldsymbol v_0\large|_{\Ga_0}=0\qquad t\in(0,+\infty),\cr
 \boldsymbol \sigma(\boldsymbol v_0)\nu\large|_{\Ga_1}=0\qquad t\in(0,+\infty),\cr \boldsymbol v_{0t}=0\qquad(x,t)\in \omega \times
(0,T),\hspace{2.8cm} (x,t)\in \omega\times(0,T).\end{cases}
\end{equation}
Then it follows from Proposition \ref{elastic-f.4} that
\be \boldsymbol v_0=0,\quad x\in \Om, \ee
which contradicts (\ref{wg-kl.89}).$\Box$

{\bf Proof of Theorem \ref{t1.2}}
It follows from  (\ref{wg-kl.23}) and (\ref{wg-kl.30}) that, for sufficiently large $T$,
\begin{eqnarray}
  TE(T)\leq \int_0^T E(t)dt&& \leq C(E_0,T)\int_0^T \int_{\Om}a(x)|\boldsymbol u_{t}|^2dx dt\nonumber\\
&& \leq C(E_0,T)(E(0)-E(T)).
 \end{eqnarray}
Hence,
\be E(T)\leq \frac{C(E_0,T)}{C(E_0,T)+T}E(0).\ee
The estimate (\ref{ex.1}) holds.$\Box$

\end{document}